\documentclass[11pt]{article}

\usepackage{fullpage}
\usepackage{graphicx} 
\usepackage{amssymb,dsfont,color}
\usepackage{amsthm}
\usepackage{soul}
\usepackage{mathtools}

\usepackage{amsmath}
\usepackage{mathrsfs}
\usepackage{enumerate}
\usepackage{mathtools}
\usepackage{stmaryrd}

\usepackage{enumitem}

\def\R{\mathbb{R}}
\def\e{\varepsilon}
\def\S{\cal{S}}

\def\ds{\displaystyle}

\def\hatlam{\hat\lambda} 
\def\lam{\lambda} 

\newtheorem{theorem}{Theorem}
\newtheorem{proposition}{Proposition}
\newtheorem{assumption}{Assumption}
\newtheorem{remark}{Remark}
\newtheorem{lemma}{Lemma}
\newtheorem{corollary}{Corollary}
\newtheorem{definition}{Definition}

\begin{document}

\title{\color{black}A Berger--Wang formula for impulsive 
switched systems}

\author{Yacine~Chitour\thanks{Yacine Chitour is with Laboratoire des Signaux et Syst\`emes (L2S), Universit\'e Paris-Saclay, CNRS, CentraleSup\'elec, Gif-sur-Yvette, France,  {\tt yacine.chitour@l2s.centralesupelec.fr}.}\quad
Jamal~Daafouz\thanks{Jamal Daafouz is with CNRS, CRAN, Université de Loraine, Institut Universitaire de France (IUF), Vandoeuvre-lès-Nancy, France, {\tt Jamal.Daafouz@univ-lorraine.fr}.}\quad
Ihab~Haidar\thanks{Ihab Haidar is with Laboratoire d'\'Electronique et d'Automatique (L\'EA), Cergy, France, and with Laboratoire Jacques-Louis Lions (LJLL), Inria, Sorbonne Universit\'e, Paris, France, {\tt ihab.haidar@ensea.fr}.} \quad
Paolo~Mason\thanks{Paolo Mason \textcolor{blue}{is} with Laboratoire des Signaux et Syst\`emes (L2S), Universit\'e Paris-Saclay, CNRS, CentraleSup\'elec, Gif-sur-Yvette, France,  {\tt paolo.mason@l2s.centralesupelec.fr}.} \quad 
Mario~Sigalotti\thanks{Mario Sigalotti is with Laboratoire Jacques-Louis Lions (LJLL), Inria, Sorbonne Universit\'e, Universit\'e de Paris, CNRS, Paris, France, {\tt mario.sigalotti@inria.fr}.}
}

\maketitle

\begin{abstract}
This paper addresses a class of impulsive systems defined by a mix of continuous-time and discrete-time switched linear dynamics.  
We first analyze a related class of weighted discrete-time switched systems for which we establish 
a Berger--Wang-type 
result. 
An analogous result is then derived for impulsive systems and subsequently used to characterize their exponential stability through a spectral approach, thereby extending existing results in switched-systems theory.
\end{abstract}


\section{Introduction}

In the present paper we focus on characterizing  {\color{black}uniform stability properties} of impulsive linear  switched systems in $\R^d$ described by  
 \begin{equation}\label{impulsive}
 \left\{\begin{array}{ll}
 \dot x(t)=Z_1(t_k) x(t),   &t\in [t_{k},t_{k+1}),\vspace{5pt}\\
    x(t_{k+1})=Z_2(t_k) x(t_{k+1}^{-}),&  k\ge 0,
\end{array}\right.
\end{equation} 
where $t\mapsto (Z_1(t),Z_2(t))$ is a piecewise-constant function taking values in a bounded set ${\cal Z}$ {\color{black}of pairs of $d\times d$ matrices. 
Here $Z$ plays the role of switching signal, $t_0=0$,  and 
$(t_k)_{k\ge 1}$ is the 
strictly increasing sequence of switching times of $Z$.

The literature on impulsive and switched systems is vast and has developed significantly over the past decades, covering linear, nonlinear, and even infinite-dimensional settings. Foundational contributions establish Lyapunov-based stability criteria, dwell-time conditions, and hybrid modeling frameworks, and provide a broad theoretical background for many classes of hybrid and switched dynamics. For general background on switched and impulsive dynamics,  including classical results on stability analysis and dwell-time techniques, we refer the reader to the monographs \cite{Haddad06} and \cite{LI2005}. A substantial body of work specifically concerns the stability of systems of the form \eqref{impulsive}. For instance, \cite{MR2405105} develops a Lyapunov characterization of exponential stability in a nonlinear setting, with applications to uncertain sampled–data control. In \cite{MR2836106}, necessary and sufficient conditions for several notions of stability are obtained in the case where the switching times $(t_k)$ are fixed. The case of positive dwell time has been investigated in a series of contributions by Briat \cite{MR3569394,MR3624532}. 
Results relying stability of impulsive switched systems and ISS properties can be found in \cite{ahmed2024lyapunov,mancilla2020uniform,Li18}. 
Further developments and more recent results can be found in \cite{MR4446170,Ren21,MR4794983}. 
Our purpose is to complement this broad context
of existing approaches by developing a spectral framework for impulsive linear switched systems.
In particular, we establish a Berger--Wang-type result adapted to systems of the form \eqref{impulsive} and show how it leads to a spectral characterization of exponential stability.
}

In order to characterize the stability of systems of the form~\eqref{impulsive}, we first examine a general class of discrete-time switched systems, referred to as {\it weighted discrete-time switched systems}. More precisely, given a {\color{black}set ${\cal N}$
made of pairs $(N,\tau)$ where $N$ is a $d\times d$ matrix and $\tau$ a nonnegative scalar, we consider the following class of systems
\begin{equation}\label{weighted}
x(k+1) = N_k x(k), \,\, (N_k,\tau_k) \in {\cal N}, \, k\ge 0,
\end{equation}
where the transition from $x(k)$ to $x(k+1)$ takes a time duration $\tau_k$, also called \emph{weight of the $k$th mode}}.
A standard discrete-time switched system can be seen as a special case of system~\eqref{weighted}, in which each mode has a unit weight. The notion of exponential stability of weighted discrete-time switched systems is defined in a manner analogous to the classical (unit-weight) case. It is important to note that the stability and instability of systems such as~\eqref{weighted} are independent of the weights (cf.~\cite{chitour2021switching}). However, the exponential growth rate, which is formally defined and studied in~\cite{chitour2021switching}, depends fundamentally on the associated weights. {\color{black}
For this class of systems, we establish a Berger–Wang–type formula \cite{berger1992bounded} which, under a suitable irreducibility assumption, equates two notions of asymptotic stability for system \eqref{weighted}: one based on the operator norm and the other on the spectral radius.}
This result plays a key role in the stability analysis of system~\eqref{weighted}. Similar types of results have been studied in \cite{WIRTH200217} in the context of discrete inclusions, for linear switched dynamical systems on graphs in~\cite{CICONE2018165}, and have been extended in \cite{dai2014robust,kozyakin2014berger} to Markovian systems. 
{\color{black}
Berger--Wang type formulas have both theoretical and computational interest: 
they provide a  characterization of the asymptotic exponential rate 
based on periodic signals 
that 
allows a structural understanding of growth phenomena
and  justifies 
numerical methods based on  spectral techniques \cite{jungers,jungers-jsr-toolbox}. 
An important theoretical application of the Berger--Wang formula is the continuity of the maximal Lyapunov exponent, as presented in Section~\ref{ssec:continuity}.  
Another motivation for us comes from~\cite{singular-preprint}, where stability criteria for hybrid linear systems subject to singular perturbations are obtained thanks to the stability analysis of auxiliary systems of the form~\eqref{impulsive}.
}

Returning to system~\eqref{impulsive}, we consider its associated weighted discrete-time switched system with modes from 
\[{\cal N}=\{(Z_2e^{tZ_1},t)\mid t\geq 0, (Z_1,Z_2)\in {\cal Z}\}.\] 

{\color{black}
Using the Berger–Wang formula developed for weighted discrete-time switched systems, and extending it to the reducible case, we characterize the stability of system \eqref{impulsive} via the sign of its maximal Lyapunov exponent: system \eqref{impulsive} is exponentially stable if and only if this exponent is negative, and exponentially unstable if and only if it is positive.

The paper is organized into four main sections. Section 2 introduces the class of impulsive linear switched systems under consideration and states the main theoretical results concerning their maximal Lyapunov exponent and stability. Section 3 studies general weighted discrete-time switched systems and establishes a Berger–Wang–type formula for their maximal Lyapunov exponent, results that form the technical core of the paper.
Section 4 focuses on impulsive switched systems: we first demonstrate how these systems can be naturally represented within the weighted discrete-time framework, and then apply and extend the results of Section 3 to derive a spectral characterization of exponential stability. The paper ends with a conclusion.
}

\subsection{\bf Notation:}\label{sec:notation}
By $\R$ we denote the set of real numbers and by $\R_{\geq \tau}$ the set of real numbers greater than $\tau\geq 0$. We use $\mathbb{N}$ for the set of positive integers. 
We use $M_{n,m}(\R)$ to denote the set of $n\times m$ real matrices and simply $M_n(\R)$ if $n=m$. The $n\times n$ identity matrix is denoted by $I_n$. 
The spectral radius of a square matrix $M$ (i.e., the maximal modulus of its eigenvalues) is denoted by $\rho(M)$ and its spectral abscissa (i.e., the maximal real part of its eigenvalues) by $\alpha(M)$. The Euclidean norm of a vector $x\in \R^n$ is denoted by $|x|$, while $\|\cdot\|$ denotes the induced norm on $M_n(\R)$, that is, $\|M\|=\max_{x\in \R^n\setminus\{0\}}\frac{|Mx|}{|x|}$ for $M\in M_n(\R)$. 
The vector space generated by a set of vectors $S$ is denoted by 
$\mathrm{span}(S)$, and $\mathrm{ker}(A)$ denotes the kernel of a matrix $A$.

Given a set ${\cal Z}$, we denote by $\S_{{\cal Z}}$ the set of right-continuous piecewise-constant functions from $\R_{\geq 0}$ to ${\cal Z}$, that is, 
those functions $Z:\R_{\geq 0}\to {\cal Z}$ such that
there exists 
an increasing sequence $(t_k=t_k(Z))_{k\in \Theta^\star(Z)}$ of switching times in $(0,+\infty)$
which is 
locally finite (i.e., has no finite density point)  and for which $Z|_{[t_k,t_{k+1})}$ is constant 
for $k,k+1\in \Theta^\star(Z)$ (with $Z|_{[0,t_1)}$ and $Z|_{(\sup_{k\in \Theta^\star(Z)}t_k,+\infty)}$ also constant).  Here 
$\Theta^\star(Z)=\emptyset$, 
 $\Theta^\star(Z)=\{1,\dots,N\}$, 
or $\Theta^\star(Z)=\mathbb{N}$, 
depending on whether $Z$ has no, $N\in \mathbb{N}$, or infinitely many switchings, respectively.  Set $t_0=0$ and, when $\Theta^\star(Z)$ is finite with cardinality $N$, 
$t_{N+1}=+\infty$.  Notice that it is allowed that the value of $Z$ is the same on two subsequent intervals between switching times.

Given $\tau\ge 0$, we denote by $\S_{{\cal Z},\tau}\subset \S_{{\cal Z},0}=\S_{\cal Z}$ the set of  piecewise-constant signals  with dwell time $\tau\ge 0$ (i.e., such that $t_{k+1}\ge t_k+\tau$ for $k
\in \Theta(Z):=\{0\}\cup \Theta^\star(Z)$). The Hausdorff distance between two nonempty subsets $X$ and $Y$ of $\R^n$ is the quantity defined by 
 \[d_H(X,Y)=\max\left\{\sup_{x\in X}d(x,Y), \sup_{y\in Y}d(y,X)\right\},\] 
where $d(x,Y)=\ds\inf_{y\in Y}|x-y|$ and $d(y,X)=\ds\inf_{x\in X}|x-y|$.

\section{Problem statement and main results}\label{sec:statement}
Let $d\in \mathbb{N}$ and ${\cal Z}$ be a bounded subset of  $M_d(\R)\times M_d(\R)$.
Consider the linear switched systems with state jumps 
\begin{equation}\label{general}
    \Sigma_{{\cal Z},\tau}:
    \left\{\begin{aligned}
 \dot x(t)&=Z_1(t_k)x(t), \quad t\in [t_k,t_{k+1}), k\in  \Theta(Z),\vspace{5pt}\\
    x(t_k)&= Z_2(t_{k-1})\ds\lim_{t\nearrow t_k}x(t),\quad k\in \Theta^\star(Z),  
\end{aligned}\right.
\end{equation} 
for $Z\in \S_{{\cal Z},\tau}$, where $\Theta(Z)$ and $\Theta^\star(Z)$, introduced in Section~\ref{sec:notation}, 
are used to parameterize the switching instants of the {\color{black}\emph{switching signal}} $Z$. We denote by $\Phi_{Z}(t,0)$ the flow from time $0$ to time $t$ of $\Sigma_{{\cal Z},\tau}$ associated with the switching signal $Z$, i.e., the matrix such that  $x_0\mapsto\Phi_{Z}(t,0)x_0$ maps the initial condition $x(0)=x_0$  to the 
evolution at time $t$ of the
corresponding solution of  $\Sigma_{{\cal Z},\tau}$. 

\begin{definition}\label{0-GES def} 
System $\Sigma_{{\cal Z},\tau}$ is said to be
\begin{enumerate}
\item
\emph{exponentially  stable} (ES, for short) if there exist $c,\delta>0$ such that \begin{equation}\label{eq-ES}
\|\Phi_{Z}(t,0)\|\leq c e^{-\delta t},\quad \forall\,t\geq 0, \forall\, Z\in {\cal S}_{{\cal Z},\tau};
\end{equation}

\item
\emph{exponentially  unstable} (EU, for short) if there exist $c,\delta>0$, $Z\in {\cal S}_{{\cal Z},\tau}$, and $x_0\in\R^d\backslash\{0\}$ such that
\[|\Phi_{Z}(t,0)x_0|\ge c e^{\delta t}|x_0|,\qquad \forall t\ge 0.\]
\end{enumerate}
The \emph{maximal Lyapunov exponent} of $\Sigma_{{\cal Z},\tau}$ is defined as
\begin{equation*}
    \lambda(\Sigma_{{\cal Z},\tau})=\limsup_{t\to+\infty}\frac{1}{t}\sup_{Z\in {\cal S}_{{\cal Z},\tau}}\log(\|\Phi_{Z}(t,0)\|),
\end{equation*}
with the convention that $\log(0)=-\infty$.
\end{definition}

We prove the following results.

\begin{theorem}\label{thm:supsup}
Let $\lambda(\Sigma_{{\cal Z},\tau})<+\infty$. Then 
{\color{black}\begin{align*}
    \lambda(\Sigma_{{\cal Z},\tau})=&\max\Big\{\sup_{(Z_1,Z_2)\in {\cal Z}}\alpha(Z_1),\\
&\sup_{Z\in {\cal S}_{{\cal Z},\tau},\;k\in \Theta^\star(Z)}\frac{\log(\rho(\Phi_Z(t_k,0)))}{t_k}\Big\}. 
\end{align*}
}
\end{theorem}

{\color{black}\begin{corollary}\label{thm:ESEU} 
Let $\lambda(\Sigma_{{\cal Z},\tau})<+\infty$. Then 
\begin{enumerate}
    \item \label{stability} System $\Sigma_{{\cal Z},\tau}$ is ES if and only if $ \lambda(\Sigma_{{\cal Z},\tau})<0$;
    \item \label{instability} System $\Sigma_{{\cal Z},\tau}$ is EU if and only if $ \lambda(\Sigma_{{\cal Z},\tau})>0$.
\end{enumerate} 
\end{corollary}

It is an open question whether one can remove the assumption 
$\lambda(\Sigma_{{\cal Z},\tau})<+\infty$ in Corollary~\ref{thm:ESEU}, that is, if $\lambda(\Sigma_{{\cal Z},\tau})=+\infty$ 
implies that System $\Sigma_{{\cal Z},\tau}$ is EU. 
In the case of linear switched systems without impulses, the boundedness of $\mathcal{Z}$ clearly implies that  $\lambda(\Sigma_{{\cal Z},\tau})<+\infty$, and the conclusion of Corollary~\ref{thm:ESEU} is well known (see, for instance, \cite[Chapter~1]{Thebook}).

}

\section{\color{black}Maximal Lyapunov exponents for weighted discrete-time switched systems}
\label{sec:wighted}

{\color{black}
In this section we present some technical results concerning  weighted discrete-time systems. These results 
 will be applied in Section~\ref{sec: impulsive systems} to impulsive switched systems, after a suitable representation of the dynamics of the latter in terms of 
weighted discrete-time systems.

}

Let ${\cal N}$ be a subset of $M_d(\R)\times [0,+\infty)$. 
We denote by $\Omega=\Omega_{\mathcal{N}}$  the set of all sequences $\omega=((A_n,\tau_n))_{n\in\mathbb{N}}$ in ${\cal N}$ such that $\sum_{k\in \mathbb N}\tau_k=+\infty$.  
For every $\omega\in \Omega$ and $k\in \mathbb N$, let 
$\omega_k$ be the finite sequence made of the first $k$ elements of $\omega$ and, using the notation for $\omega$ introduced above,  associate with $\omega_k$ the weight $|\omega_k|=\tau_1+\dots+\tau_k$
and the matrix product $\Pi_{\omega_k}=A_k\cdots A_1$.
For $k_1\le k_2$, we define also the matrix product $\Pi_{\omega_{k_1\to k_2}}=A_{k_2}\cdots A_{k_1+1}$,
with the convention that $\Pi_{\omega_{k_1\to k_2}}=I_d$  if $k_1=k_2$.

We associate with $\mathcal{N}$ a system $\Xi=\Xi_{\mathcal{N}}$ whose
trajectories are the sequences $(x(k))_{k\in \mathbb N}$ in $\mathbb R^d$
such that there exists $\omega\in \Omega$ for which 
\begin{equation*}
x(k)=\Pi_{\omega_k}x(0),\qquad k\in \mathbb N.
\end{equation*}
We say that  $\Xi$ is a \emph{weighted discrete-time system}. 
{\color{black}Such systems have already been introduced in \cite{chitour2021switching} with the restriction that all the weights are strictly positive. In the present paper, we extend the framework to nonnegative weights. }

\begin{remark}
As we will see in Section~\ref{sec: impulsive systems}, we can associate with  a switched system with jumps a 
natural weighted discrete-time system. 
Weighted discrete-time systems can be used to study also more general classes of systems: for instance, 
we could consider the case where the dwell time depends on the mode and has also an upper bound, in the sense that there exist
$\tau_-,\tau_+:{\cal Z}\to \R_{\ge 0}\cup\{+\infty\}$ with $\tau_-\le \tau_+$ such that $t_k+\tau_+(Z(t_k))\ge t_{k+1}\ge t_k+\tau_-(Z(t_k))$ 
for every $k\in\Theta^\star(Z)$. (See, for instance, \cite{dellarossa,protasov-kamalov}). 
\end{remark}

\begin{definition}\label{def:discret_weight}
We say that $\Xi$ is 
\begin{enumerate}
\item
\emph{exponentially  stable} (ES, for short) if there exist $c,\delta>0$ such that
\begin{equation}\label{eq-ES-discret}
\|\Pi_{\omega_k}\|\leq c e^{-\delta |\omega_k|},\qquad \forall\,\omega\in\Omega,\; \forall\,k\in \mathbb{N};
\end{equation}

\item
\emph{exponentially  unstable} (EU, for short) if there exist $c,\delta>0$, $x_0\in\R^d\setminus\{0\}$, and $\omega\in \Omega$   such that 
\[|\Pi_{\omega_k}x_0|\ge c e^{\delta |\omega_k|}|x_0|, \quad \forall\,k\in \mathbb{N}.\]
\end{enumerate}

The \emph{maximal Lyapunov exponent} of a weighted discrete-time switched system $\Xi=\Xi_{\mathcal{N}}$ is defined as

\begin{equation}\label{Lyap_exp_graph}
\lam(\Xi)=\limsup_{t\to+\infty}\sup_{
\omega\in \Omega,\, k\in \mathbb N, \,  
|\omega_k|=t}\dfrac{\log(\|\Pi_{\omega_k}\|)}{|\omega_k|},
\end{equation}
with the convention that $\sup\emptyset=-\infty$ and $\log 0 = -\infty$.
\end{definition}

We define also
\begin{equation}\label{Lyap_exp_mode}
\hatlam(\Xi)=\limsup_{t\to+\infty}\sup_{
\{A\mid (A,t)\in {\cal N}\}}\dfrac{\log(\|A\|)}{t}.
\end{equation}

We now define the {\color{black}\emph{generalized spectral abscissa}} $\mu(\Xi)$, which is the counterpart of $\lam(\Xi)$ in which the norm of a product $\Pi_{\omega_k}$ is replaced by its spectral radius, that is,
\begin{equation}\label{eq:tildemu-def}
\mu(\Xi)=\limsup_{t \to +\infty} \sup_{\omega \in \Omega,\, k \in \mathbb{N},\, |\omega_k| = t} \frac{\log(\rho(\Pi_{\omega_k}))}{|\omega_k|}.
\end{equation}

\begin{remark}\label{rmk:rho-sup-limsup}
It holds that $\mu(\Xi)
=\sup_{\omega\in \Omega,\;k\in\mathbb{N}}  \frac{\log(\rho(\Pi_{\omega_k}))}{|\omega_k|}$.
Indeed, by definition of $\mu(\Xi)$,  one immediately has
$\mu(\Xi)
\le \sup_{\omega\in \Omega,\;k\in\mathbb{N}}  \frac{\log(\rho(\Pi_{\omega_k}))}{|\omega_k|}$.
Conversely, 
fix $\omega\in \Omega$ and $k\in  \mathbb{N}$, and denote by $\omega_k^{*}\in \Omega$ the repetition of the finite sequence $\omega_k$ for infinitely many times. Then, for every $m\in \mathbb{N}$, $\rho(\Pi_{(\omega_k^{*})_{mk}})=\rho(\Pi_{\omega_k}^{m})=\rho(\Pi_{\omega_k})^m$
and $|(\omega_k^{*})_{mk}|=m |\omega_k|$, whence the inequality $\frac{\log(\rho(\Pi_{\omega_k}))}{|\omega_k|}\le \mu(\Xi)$. 
\end{remark}

Let, for $\xi\in \R$, $\Omega^{\xi}$ be the set of sequences $((e^{\xi\tau_n}A_n,\tau_n))_{n\in \mathbb{N}}$ such that $(A_n,\tau_n)\in {\cal N}$ and denote by $\Xi^\xi$ the corresponding discrete-time weighted system.

\begin{lemma}\label{shift}
For every $\xi\in \R$ we have  
\begin{equation*}
\lam(\Xi^\xi)=\xi+\lam(\Xi),\; \hatlam(\Xi^\xi)=\xi+\hatlam(\Xi),\; \mu(\Xi^\xi)=\xi+\mu(\Xi). 
\end{equation*}
\end{lemma}

\begin{proof}
The proof is direct from the definition of $\lam(\Xi)$, $\hatlam(\Xi)$ and $\mu(\Xi)$.
\end{proof}

Let us recall the definition of 
irreducible set of matrices. 
\begin{definition}\label{invariant}
We say that a subspace $E$ of $\R^d$ is invariant with respect to a set ${\cal M}\subset M_{d}(\R)$ if and only if it is invariant with respect to every matrix $M\in {\cal M}$, i.e.,  $Mx\in E$ for every $x\in E$. A set ${\cal M}\subset M_{d}(\R)$ is said to be irreducible if its only invariant subspaces are $\{0\}$ and $\R^d$. Otherwise it is said to be reducible.
\end{definition}

We will make use of the following result.
\begin{lemma}\label{lemma-positivetime}
Let $\mathcal{M} \subset M_d(\mathbb{R})$ be irreducible and consider a nonzero matrix $\bar M\in \mathcal{M}$. Then, for every $x\in \mathbb{R}^d\setminus\{0\}$ there exists a product $\Pi$ of matrices in $\mathcal{M}$ containing $\bar M$ among its factors such that $\Pi x \neq 0$. 
\end{lemma}
\begin{proof}
In the following we denote as $\mathcal{M}^k$ the set of all possible products of $k$ matrices of $\mathcal{M}$ and by $\mathcal{M}^k x$ the set of all possible evaluations of matrices of $\mathcal{M}^k$ at $x\in \mathbb{R}^d$.

Let $x\in \mathbb{R}^d\setminus\{0\}$.
 If $x\notin  \mathrm{ker}(\bar M)$ then we can take $\Pi = \bar M$.
Assume now $x\in \mathrm{ker}(\bar M)$. By irreducibility of $\mathcal{M}$ we have that the vector space $\mathrm{span}\{\mathcal{M}^kx\mid k\in\mathbb{N}\}$, which is invariant with respect to $\mathcal{M}$, coincides with $\mathbb{R}^d$. This means that there exist $k\in\mathbb{N}$ and $\bar \Pi\in \mathcal{M}^k$ such that $\bar \Pi x\notin  \mathrm{ker}(\bar M)$ (note that $ \mathrm{ker}(\bar M)\neq \mathbb{R}^d$ since $\bar M$ is nonzero).
We conclude by taking 
$\Pi = \bar M\bar \Pi$.
\end{proof}

For general weighted discrete-time switched systems, we need the following assumption.

\begin{assumption}\label{irreducible}
The {\color{black}projection onto $M_d(\R)$ of the set ${\cal N}$, i.e,} 
$\{A\mid (A,\tau)\in {\cal N}
\},$
is irreducible. 
\end{assumption}

{\color{black}
The following proposition states that, in the irreducible case, the Lyapunov exponent 
$\lam(\Xi)$ provides a uniform upper bound of the exponential growth 
of the trajectories of $\Xi$.

}

\begin{proposition}\label{quasi-barab-0}
Assume that $\lam(\Xi)$  is finite
and that Assumption~\ref{irreducible} holds. 
Then there exists $C\geq 1$ such that  
    $\|\Pi_{\omega_k}\|
    \leq C e^{\lam(\Xi)|\omega_k|}$ for every $\omega\in {\Omega}$ and $k\in \mathbb{N}$.
\end{proposition}
\begin{proof}
Since $\lam(\Xi)$ is finite, we can suppose without loss of generality that $\lam(\Xi)=0$ by Lemma~\ref{shift}. 
Consider the set 
\[E:=\{x\in\R^d \mid  \sup_{\omega\in\Omega, k\in \mathbb{N} }|\Pi_{\omega_k}x|<\infty\}.\]
The set $E$ is a subspace of $\R^d$ which is invariant with respect to 
$\{A\in M_d(\mathbb{R})\mid (A,\tau)\in {\cal N}
\}$, which is irreducible by Assumption~\ref{irreducible}. Hence, $E$ is either $\{0\}$ or $\R^d$. 

We will prove that $E=\mathbb{R}^d$.  By contradiction, assume that $E=\{0\}$. Then we claim that there exists a finite subset $\mathcal{W}$ of $\Omega\times \mathbb{N}$ such that, for every nonzero $x\in\R^d$, there exists $(\omega,k)\in\mathcal{W}$ such that $|\omega_{k}|>0$ and $|\Pi_{\omega_{k}}x| > 2|x|$. Note that we can choose the same pair $(\omega,k)$ for all vectors proportional to $x$. Then, we just need to prove the claim for all $x$ belonging to the unit sphere of $\mathbb{R}^d$. To prove the claim we first apply Lemma~\ref{lemma-positivetime} with $\mathcal{M} = \{A\in M_d(\mathbb{R})\mid (A,\tau)\in {\cal N}\}$ and $\bar M\in  \{A\in M_d(\mathbb{R})\mid (A,\tau)\in {\cal N},\tau>0\}$ (the fact that such a nonzero $\bar M$ exists follows from the fact that, by assumption, $\lam(\Xi)>-\infty$)
and deduce 
 that, for every $x\in\R^d$ with $|x|=1$ there exists  $\omega^1\in\Omega$ and $k^1\in\mathbb{N}$ such that $|\omega^1_{k^1}|>0$ and $\Pi_{\omega^1_{k^1}}x\neq 0$. Then, since $E=\{0\}$, there exists $\omega^2\in\Omega$ and $k^2\in\mathbb{N}$ such that $|\Pi_{\omega^2_{k^2}}\Pi_{\omega^1_{k^1}}x| > 2|x|$. In other words we have constructed $\omega\in \Omega$ such that, setting $k= k^1+k^2$ we have $|\Pi_{\omega_{k}}x|> 2|x|$ and $|\omega_{k}| = |\omega^1_{k^1}| + |\omega^2_{k^2}| >0$. 
For $x\in\R^d$ with $|x|=1$ we can find $(\omega^x,k^x)\in\Omega\times \mathbb{N}$ such that $|\omega^x_{k^x}|>0$ and $|\Pi_{\omega^x_{k^x}}y| > 2|y|$ for every $y\in U_{\omega^x,k^x}$, $U_{\omega^x,k^x}$ being a open neighborhood of $x$. The sets  $U_{\omega^x,k^x}$ form an open covering of the unit sphere and, by compactness of the latter, we can extract a finite covering associated with a finite set $\mathcal{W}$ of elements of $\Omega\times \mathbb{N}$ which satisfies the claim.

Let $\delta = \min_{(\omega,k)\in \mathcal{W}}|\omega_k|$ and $\Delta = \max_{(\omega,k)\in \mathcal{W}}|\omega_k|$. 
Fix $x_0\in \R^d\setminus\{0\}$. Then we can construct recursively a sequence $\bar \omega\in \Omega$ by concatenating finite sequences $\omega^i_{k^i}$, with $(\omega^i,k^i)\in \mathcal{W}$ for $i\in\mathbb{N}$, in such a way that, setting $x_i=\Pi_{\omega^{i-1}_{k^{i-1}}}\cdots \Pi_{\omega^{0}_{k^{0}}}x_{0}$, one has $|\Pi_{\omega^i_{k^i}}x_i|\ge 2|x_i|$. Setting $\ell^n = \sum_{i=0}^{n-1} |\omega^i_{k^i}|$ we then have $\lim_{n\to +\infty}  |\bar\omega_{\ell^n}|\geq \lim_{n\to +\infty}  n\delta = +\infty$ so that
\begin{align*}
0 = \lam(\Xi)&\ge\limsup_{n\to+\infty}\frac{\log\|\Pi_{\bar\omega_{\ell^n}}\|}{|\bar\omega_{\ell^n}|}\geq \limsup_{n\to+\infty}\frac{\log |\Pi_{\bar\omega_{\ell^n}}x_0|}{|\bar\omega_{\ell^n}|}\\
&\geq  \limsup_{n\to+\infty}\frac{\log (2^n |x_0|)}{n\Delta} = \frac{\log 2}{\Delta},
\end{align*}
which is a contradiction. We have therefore proved that  $E$ is equal to $\R^d$. 

Now, assume by contradiction that $\sup_{\omega\in\Omega, k\in \mathbb{N}}\|\Pi_{\omega_k}\|$ is not finite. Then, there exist a sequence $\{(\omega^n,k^n)\}_{n\in \mathbb{N}}$ of elements of $\Omega\times\mathbb{N}$ and a sequence of unit vectors $\{x^n\}_{n\in \mathbb{N}}$ in $\R^d$ so that 
$|\Pi_{\omega^n_{k^n}}x^n|=\|\Pi_{\omega^n_{k^n}}\|$ tends to infinity as $n$ goes to infinity. Up to a subsequence, one has that $\lim_{n\to +\infty}x^n=x^*\in\R^d$, and 
\begin{align*}
|\Pi_{\omega^n_{k_n}}x^*|&\ge |\Pi_{\omega^n_{k^n}}x^n|-|\Pi_{\omega^n_{k^n}}(x^n-x^*)|\\
&\ge \|\Pi_{\omega^n_{k_n}}\|\left(1-|x^n-x^*|\right)\ge \frac{\|\Pi_{\omega^n_{k^n}}\|}{2}, 
\end{align*}
where the last inequality holds true for $n$ large enough. Hence, $x\notin E$, which is a contradiction with the fact that $E=\R^d$. Hence $\sup_{\omega\in\Omega, k\in \mathbb{N}}\|\Pi_{\omega_k}\|$ is finite, which concludes the proof of the proposition.
\end{proof}

{\color{black}
In the next proposition we prove that, under the additional assumption $\hatlam(\Xi)<\lam(\Xi)$,
the upper bound $\lambda(\Xi)$ on the uniform exponential growth rate of the trajectories of $\Xi$ is sharp in a suitable sense. 

}

\begin{proposition}\label{quasi-barab-1}
Assume that $\hatlam(\Xi)<\lam(\Xi)<+\infty$
and that Assumption~\ref{irreducible} holds. 
Then for every $t>0$ and every $x\in\mathbb{R}^n$, there exists $(\omega,k)\in\Omega\times\mathbb{N}$ satisfying $|\omega_k|\geq t$ and $|\Pi_{\omega_k}x| \ge c e^{\lam(\Xi)|\omega_k|}|x|$ for some positive constant $c$ only depending on $\mathcal{N}$.
\end{proposition}
\begin{proof}
As in the proof of Proposition~\ref{quasi-barab-0} we assume, without loss of generality, that $\lam(\Xi)=0$.
Consider the set 
\[
{\cal R}_{\infty}:=\cap_{t\geq 0}{\cal R}_{t},
\]
where $\mathcal{R}_t  = \overline{\left\{\Pi_{\omega_{k}}\mid (\omega,k)\in\Omega\times \mathbb{N}\mbox{ s.t. } |\omega_{k}|\ge t\right\}}.$
The sequence of sets $\mathcal{R}_t$ is decreasing, in the sense that $\mathcal{R}_t\supset \mathcal{R}_s$ whenever $s\geq t$, and each $\mathcal{R}_t$ is 
 closed and, by Proposition~\ref{quasi-barab-0}, bounded. Hence ${\cal R}_{\infty}$ is closed and bounded, that is, it is a compact subset of $M_d(\mathbb{R})$. Moreover, ${\cal R}_{\infty}$ is nonempty by Cantor intersection theorem.

We claim that  ${\cal R}_{\infty}\ne \{0\}$. 
Assume by contradiction that ${\cal R}_{\infty}=\{0\}$. Then 
\begin{align*}
\cap_{t\geq 0}\left( \mathcal{R}_t \cap \{A\in M_d(\mathbb{R})\mid \|A\|\ge 1/2\}\right) =\\
 {\cal R}_{\infty} \cap \{A\in M_d(\mathbb{R})\mid \|A\|\ge 1/2\} = \emptyset.
\end{align*}
It follows from Cantor intersection theorem that $\mathcal{R}_T \cap \{A\in M_d(\mathbb{R})\mid \|A\|\ge \frac12\}= \emptyset$ for some $T$ large enough, that is,
\[
    \|\Pi_{\omega_k}\|<\frac{1}2,\qquad \forall\,(\omega,k)\in \Omega\times\mathbb{N}\mbox{ s.t. }|\omega_k|\ge T. 
    \]
Next, we derive further estimates of $\|\Pi_{\omega_k}\|$ for finite sequences $\omega_k$ in two special cases.

$(A)$ Consider a finite sequence $\omega_k=\{(A_n,\tau_n)\}_{n=1,\dots,k}$ with $\tau_n < T$ for every $n$. Then  $\omega_k$ can be written as the concatenation of  finite sequences $\omega^i_{k^i}$, $i=1,\dots,\ell$ (for some $\ell\geq 0$), of elements of $\mathcal{N}$ such that $|\omega^i_{k^i}|\in [T,2T)$ together with a sequence $\bar \omega_{\bar k}$  such that $|\bar \omega_{\bar k}|\leq T$. In particular $\ell$ satisfies $|\omega_k|/(2T)-1/2\leq \ell \leq |\omega_k|/T$.
Then 
\begin{equation}\label{estimate1}
\|\Pi_{\omega_k}\| \leq \|\Pi_{\bar\omega_{\bar k}}\|\prod_{i=1}^{\ell}\|\Pi_{\omega^i_{k^i}}\|\leq C2^{-\ell}\leq \sqrt{2}C2^{-\frac{|\omega_k|}{2T}},
\end{equation}
where $C$ is as in Proposition~\ref{quasi-barab-0}.

$(B)$ 
 Pick $\gamma$ in the open interval $(\hatlam(\Xi),0)$. By definition of $\hatlam(\Xi)$
and the inequality $\hatlam(\Xi)<\gamma$, 
it follows that, up to increasing
$T$, 
$\|A\|\le \frac{1}{\sqrt{2}C} e^{\gamma \tau}$ for every $(A,\tau)\in {\cal N}$ with $\tau\ge T$.
{\color{black} In particular, given a finite sequence $\omega_k=\{(A_n,\tau_n)\}_{n=1,\dots,k}$ with $\tau_n \geq T$ for every $n$, we have}
\begin{equation}\label{estimate2}
\|\Pi_{\omega_k}\|\le \frac{1}{\sqrt{2}C} e^{\gamma |\omega_k|}.
\end{equation}

We observe now that for every $\omega$ in $\Omega$ with $\lim_{k\to +\infty}|\omega_k| = +\infty$ and every positive integer $\ell$ the finite sequence $\omega_{\ell}$ can be written as a concatenation of finite sequences which alternate between type $(A)$ and type $(B)$. By submultiplicativity of the matrix norm and the estimates~\eqref{estimate1} and \eqref{estimate2} it follows that $\|\Pi_{\omega_{\ell}}\|\le \sqrt{2}C e^{-\alpha |\omega_{\ell}|}$, where $\alpha = \min \{-\gamma,\log 2/(2T)\}
>0
$. 
Since $\omega$ and $\ell$ have been chosen arbitrarily, using the definition \eqref{Lyap_exp_graph} of $\lam(\Xi)$ we obtain $\lam(\Xi)\leq-\alpha<0$, which is a contradiction.
This concludes the proof that ${\cal R}_{\infty}\ne \{0\}$.

To conclude the proof it is enough to show that the map 
 \[ 
 x\mapsto v_t(x):=\max_{R\in {\cal R}_{t}}|Rx|
\]
 is a norm on $\R^d$ for every $t\geq 0$. Indeed, in this case, by the equivalence of norms on finite-dimensional spaces, there exists $\kappa>0$ such that with every nonzero $x\in \mathbb{R}^d$ and $t\geq 0$ one can associate an element $R\in\mathcal{R}_t$ satisfying $|Rx| \geq \kappa |x|$. Then, from the definition of $\mathcal{R}_t$, it follows that there exists $(\omega,k)\in\Omega\times \mathbb{N}$ with  $|\omega_{k}|\ge t$ such that $|\Pi_{\omega_{k}}x|\geq \frac{\kappa}{2}|x|$.

Let us prove that $v_t$ is a norm. 
By compactness of $\mathcal{R}_t$ one has that $v_t$ is well defined. Furthermore, by definition, it is clear that $v_t$ is absolutely homogeneous {\color{black}(i.e., $v_t(\xi x)=|\xi|v_t(x)$ for every $x\in \R^d$ and $\xi\in \R$)} and satisfies the triangle inequality. It remains to show that $v_t$ is strictly positive outside the origin.
 Note that, by definition, each $\mathcal{R}_t$ is invariant with respect to right multiplication by elements of $\{A\in M_d(\mathbb{R})\mid (A,\tau)\in \mathcal{N}\}$. It follows that the vector space 
 \[\{x\in  \mathbb{R}^d\mid Rx=0\quad \forall R\in \mathcal{R}_t\},\]
which is a strict subspace of $\mathbb{R}^d$ since $ \mathcal{R}_t\neq \{0\}$, is invariant with respect to the set $\{A\in M_d(\mathbb{R})\mid (A,\tau)\in \mathcal{N}\}$ and therefore, by Assumption~\ref{irreducible}, it coincides with $\{0\}$. This means that $v_t$ is strictly positive outside the origin. We have therefore shown that  $v_t$ is a norm.
 
 This concludes the proof of the proposition.
\end{proof}

\begin{remark}
Proposition~\ref{quasi-barab-1} cannot be extended, in general, to  the case where $\hatlam(\Xi)=\lam(\Xi)$.
Consider, for example,  the case where 
${\cal N}=
\{(\frac1n e^{-\frac \tau n},\tau)\mid n\ge 2,\;\tau\ge 0\}$.
In this case it is easy to verify that $\hatlam(\Xi)=\lam(\Xi)=0$. 

On the other hand, consider a general $\omega\in \Omega$, that is a sequence $((\frac1{n_k} e^{-\frac {\tau_k}{n_k}},\tau_k))_{k\in\mathbb{N}}$ in $\mathcal{N}$.
Notice that 
\[\|\Pi_{\omega_k}\|\le \frac{1}{\max\{n_1,\dots,n_k\}}e^{-\frac{|\omega_k|}{\max\{n_1,\dots,n_k\}}},\qquad \forall k\in\mathbb{N}.\] 

Assume by contradiction that there exists $c>0$ as in 
Proposition~\ref{quasi-barab-1}.
In order to have
$|\Pi_{\omega_k}x|\ge c|x|$ for a given $x\ne 0$, one must have 
$\frac{1}{\max\{n_1,\dots,n_k\}}>c$, that is,  $\max\{n_1,\dots,n_k\}<\frac 1 c$. This implies that $e^{-\frac{|\omega_k|}{\max\{n_1,\dots,n_k\}}}$ is smaller than $e^{-c|\omega_k|}$, 
 so that
$e^{-c|\omega_k|}|x|\ge c|x|$. 
Since in 
Proposition~\ref{quasi-barab-1} 
the sequence $\omega$ and the integer $k$ are taken in such a way that $|\omega_k|>t$ and $t$ can be arbitrarily large, a contradiction is reached.
\end{remark}

{\color{black}
We are now ready to prove a key result for obtaining a version of the Berger--Wang formula for weighted discrete-time systems.  
}

\begin{proposition}\label{lambdatilde-rho}
Let Assumption~\ref{irreducible} hold. 
Assume, moreover, that 
$0<\lam(\Xi)<+\infty$. Then, 
either $\hatlam(\Xi)=\lam(\Xi)$ 
or there exists 
$\omega\in\Omega$ and $k\in \mathbb{N}$ such that 
$\rho(\Pi_{\omega_k})>1$.
\end{proposition}
\begin{proof}
By definitions of $\hatlam(\Xi)$ and $\lam(\Xi)$, we have either $\hatlam(\Xi)=\lam(\Xi)$ or $\hatlam(\Xi)<\lam(\Xi)$.

Assume $\hatlam(\Xi)<\lam(\Xi)$. In this case, by Lemma~\ref{shift}, we can equivalently prove that $\sup_{(\omega,k)\in \Omega\times\mathbb{N}} \rho(\Pi_{\omega_k})\ge 1$ whenever $\lam(\Xi)= 0$.
By  
Proposition~\ref{quasi-barab-0},
there exists $C>0$ such that  
\begin{equation}\label{eq:exp--bound}
\|\Pi_{\omega_k}\|\le C,\qquad \mbox{for every $\omega\in\Omega$ and $k\in\mathbb{N}$}.
\end{equation}
Moreover, by 
Proposition~\ref{quasi-barab-1}, 
there exists a sequence $(\omega^n,k_n)_n$ in $\Omega\times \mathbb{N}$ with $|\omega^{n}_{k_n}|\to +\infty$ 
such that $\|\Pi_{\omega^n_{k_n}}\|\geq c$ for some $c>0$ independent of $n$. 

Denote by $\tau_n$ the maximal weight of an element of $\omega_{k_n}^n$ and by $(A_n,\tau_n)$ the corresponding element of ${\cal N}$. We claim that  
\begin{equation}\label{eq:tn-bounded}
    \sup_n \tau_n<+\infty.
\end{equation} 
Indeed, if this were not the case,
write $\Pi_{\omega_{k_n}^n}$ as 
$\Pi_{\nu_n}A_n\Pi_{\mu_n}$, for some finite sequences $\nu_n$ and $\mu_n$ in ${\cal N}$. 
Applying \eqref{eq:exp--bound} to  $\mu_n$ and $\nu_n$ and 
using the relation $\tau_n=|\omega_{k_n}^n|-|\nu_n|-|\mu_n|$,
we deduce that 
\[\|A_n\|\ge \frac{c}{C^2}.
\]
Since $\tau_n\to +\infty$ by the contradiction assumption, we have, by definition of $\hatlam(\Xi)$, 
that
$\hatlam(\Xi)\ge 0$, which is impossible given that $\hatlam(\Xi)<\lam(\Xi)=0$. This concludes the proof of \eqref{eq:tn-bounded}. 

Let us now define the vectors 
$y_{k,n}=\Pi_{\omega^n_{k}}x_n$ for $k=1,\dots, k_n,$
where $(x_n)_n$ is a sequence of unit vectors such that  $\|\Pi_{\omega^n_{k_n}}\|= |\Pi_{\omega^n_{k_n}}x_n|$. 
Since 
\[|\Pi_{\omega_k^n\to \omega_{k_n}^n}y_{k,n}|=|\Pi_{\omega_{k_n}^n}x_n|\ge c,\]
we deduce from \eqref{eq:exp--bound} that each 
$y_{k,n}$ belongs to $K$, where 
\[K:=\left\{x\in \mathbb{R}^d\mid \frac{c}{C}\le |x|\le C\right\}.\] 

We define the set 
\begin{equation*}
I_n=\left\{(j_1,j_2)\mid 1\leq j_1<j_2\leq k_n 
 \right\}.
\end{equation*}
By \eqref{eq:tn-bounded}, 
we have $\#I_n\to +\infty$. For $n\geq 1$, let $j_1^n, j_2^n$ be such that 
\begin{equation*}
(j_1^n,j_2^n)\in {\arg\min}_{(j_1,j_2)\in I_n}|y_{j_1,n}-y_{j_2,n}|
\end{equation*} 
and notice that \begin{equation}\label{eq:lim-argmin}
\lim_{n\to +\infty}|y_{j_1^n,n}-y_{j_2^n,n}|\to 0
\end{equation}
by the boundedness of $K$ and unboundedness of $\#I_n$.
We have 
\begin{equation}\label{parallel}
y_{j_2^n,n}=\Pi_{\omega^n_{j_1^n\to j_2^n}}y_{j_1^n,n}, \quad \forall\, n\geq 1.
\end{equation}
Up to extracting a subsequence, we can assume that
$y_{j_1^n,n}$ converges to some $y^\star$. Notice that $y^\star\ne0$ (by definition of $K$) and that $\lim_{n\to\infty}y_{j_2^n,n}=y^*$ by \eqref{eq:lim-argmin}.

By \eqref{eq:exp--bound}, we can extract a  subsequence of $\Pi_{\omega^n_{j_1^n\to j_2^n}}$ 
converging to some $M\in M_d(\R)$
as $n$ tends to infinity. 
By~\eqref{parallel} we deduce that $My^\star=y^\star$, implying that $\rho(M)\ge 1$. 
Therefore,  
\begin{equation*}
\lim_{n\to\infty}\rho(\Pi_{\omega^n_{j_1^n\to j_2^n}})=\rho(M)\geq 1, 
\end{equation*}
 concluding the proof.
\end{proof}

\begin{corollary}
[Berger--Wang formula for $\Xi$]
\label{Berger-Wang}
Let Assumption~\ref{irreducible} hold and 
suppose that $\lam(\Xi)<+\infty$. 
Then  
\begin{align}
   \lam(\Xi)= \max\left\{\hatlam(\Xi),\mu(\Xi)\right\}.
\label{BW-discret}
\end{align}
\end{corollary}

\begin{proof}
Observe that $\rho(\Pi_{\omega_k}) \leq \|\Pi_{\omega_k}\|$  for every $(\omega, k) \in \Omega \times \mathbb{N}$. 
Hence, $\mu(\Xi)\le \lam(\Xi)$.
Using the inequality $\hatlam(\Xi) \le \lam(\Xi)$, we get that $ \max\{\hatlam(\Xi), \mu(\Xi)\}\le \lam(\Xi)$. 

We are left 
to show that $ \max\{\hatlam(\Xi), \mu(\Xi)\}\ge \lam(\Xi)$. 
If $\hatlam(\Xi)=\lam(\Xi)$ or $\lam(\Xi)=-\infty$ the conclusion holds true. 

Let us consider $\e>0$ and apply
Lemma~\ref{shift} with 
 $\xi=\e-\lam(\Xi)$, noticing that  $\lam(\Xi)$ is finite by assumption.
 Since $\lam(\Xi^\xi)=\e>0$, we deduce from Proposition~\ref{lambdatilde-rho}  that either $\hatlam(\Xi^\xi)=\lam(\Xi^\xi)$ or $\mu(\Xi^\xi)>0$, i.e., either $\hatlam(\Xi)=\lam(\Xi)$ or $\mu(\Xi)+\e>\lam(\Xi)$. Since this holds for any $\e>0$ and we have reduced our analysis to the case $\hatlam(\Xi)<\lam(\Xi)$, we conclude that 
 $\mu(\Xi)\ge\lam(\Xi)$. Therefore, $ \max\{\hatlam(\Xi), \mu(\Xi)\}\ge \lam(\Xi)$, which completes the proof. 
\end{proof}

{\color{black}
\section{Stability of 
impulsive systems
by means of
the weighted discrete-time formalism
}\label{sec: impulsive systems}

\subsection{System $\Sigma_{{\cal Z},\tau}$ seen as a weighted discrete-time switched system}

}

In order to study the stability of
System~$\Sigma_{{\cal Z},\tau}$ 
introduced in Section~\ref{sec:statement}, we will consider the 
weighted discrete-time switched system
$\Xi_{{\cal Z},\tau}:=\Xi_{{\cal N}_{{\cal Z},\tau}}$, where 
$
{\cal N}_{{\cal Z},\tau}:=\{(Z_2e^{tZ_1},t)\mid t\geq \tau, (Z_1,Z_2)\in {\cal Z}\}$.
Denote also $\Omega_{{\cal Z},\tau}:=\Omega_{{\cal N}_{{\cal Z},\tau}}$, where we again use the notation introduced at the beginning of Section~\ref{sec:wighted}.
It is important to note, before presenting the analysis here, that Assumption~\ref{irreducible} is no longer required for this particular class of weighted discrete-time switched systems.

The two systems $\Sigma_{{\cal Z},\tau}$ and $\Xi_{{\cal Z},\tau}$ are strongly related but not completely equivalent. 
With every $Z\in\mathcal{S}_{{\cal Z},\tau}$ with 
infinitely many switchings and every initial condition $x_0$, we can associate a 
trajectory of $\Xi_{\cal Z,\tau}$ given by the evaluation 
at the switching times
 of $Z$ of the trajectory of $\Sigma_{\mathcal{Z},
 \tau}$ starting from $x_0$ and corresponding to $Z$. 
{\color{black}However, a trajectory of $\Sigma_{\mathcal{Z},\tau}$ generated by a signal $Z \in \mathcal{S}_{\mathcal{Z},\tau}$ with finitely many switchings does not, in general, admit a corresponding trajectory of $\Xi_{\mathcal{Z},\tau}$ with the same asymptotic behavior.}
Moreover, in the case where $\tau=0$, $\Xi_{\cal Z,\tau}$ may contain more trajectories than those corresponding to trajectories of $\Sigma_{\mathcal{Z},\tau}$, since $\mathcal{N}_{\mathcal{Z},0}$ contains also elements of the type $(Z_2,0)$ for $(Z_1,Z_2)\in\mathcal{Z}$, while the distance between two switching times is always positive. (In this sense, $\mathcal{N}_{\mathcal{Z},0}$ can be used to study also switched systems that can jump several times at the same time instant, provided that there are finitely many jumps on any positive time-interval).

Notice that
\begin{equation}\label{Xi<=Sigma}
\lam(\Xi_{\mathcal{Z},\tau})\le \lambda(\Sigma_{\mathcal{Z},\tau}),
\end{equation}
as it can be deduced from the definition of the two maximal Lyapunov exponents, by noticing that  for every $\omega\in \Omega_{\mathcal{Z},\tau}$ and every $k\in \mathbb{N}$, there exist a sequence $(Z^n)_{n\in\mathbb{N}}$ in $\mathcal{S}_{{\cal Z},\tau}$ and a sequence $(t_n)_{n\in\mathbb{N}}$ in $[0,+\infty)$ such that $\lim_{n\to\infty} t_n=|\omega_k|$ and $\lim_{n\to\infty}\Phi_{Z^n}(t_n,0)=\Pi_{\omega_k}$. Each $Z^n$ can be constructed 
by associating with $\omega_k=((Z_{j,2} e^{t_j Z_{j,1}},t_j))_{j=1}^k$ the 
piecewise-constant signal 
whose $j$th piece is equal to $(Z_{j,1},Z_{j,2})$  on an interval of length $t_j+\frac1n$.

\begin{remark}\label{remark-infini}
Notice that
both $\lam(\Xi_{\cal Z,\tau})=-\infty$ and $\lam(\Xi_{\cal Z,\tau})=+\infty$ may occur. For instance, 
$\lam(\Xi_{\cal Z,\tau})=-\infty$ when $Z_2=0$ for every $(Z_1,Z_2)\in \mathcal{Z}$. 
As for $\lam(\Xi_{\cal Z,\tau})=+\infty$, a necessary condition for it to happen is that $\tau=0$. 
Indeed, if $\tau$ is positive and  $\gamma_1,\gamma_2\in (0,+\infty)$ are taken so that $\|Z_1\|\le \gamma_1$ and $\|Z_2\|\le e^{\gamma_2\tau}$ for every $(Z_1,Z_2)\in {\cal Z}$ (which is possible because $\mathcal{Z}$ is bounded), then $\|\Pi_{\omega_k}\|\le e^{(\gamma_1+\gamma_2)|\omega_k|}$ for every $k\in \mathbb{N}$ and $\omega\in \Omega_{\cal Z,\tau}$, yielding that $\lam(\Xi_{\cal Z,\tau})\le \gamma_1+\gamma_2$. 
Lemma~\ref{lem:triip} below characterizes the case where 
$\lam(\Xi_{{\cal Z},0})=+\infty$.
\end{remark}

\begin{lemma}\label{lem:triip}
The following three properties are equivalent:
\begin{enumerate}
\item \label{It1}$\lam(\Xi_{{\cal Z},0})<+\infty$;
    \item \label{It2} The set $\{Z^1_2\cdots Z^k_2\mid k\geq 1,\; (Z^1_1,Z^1_2),\dots, (Z^k_1,Z^k_2)\in {\cal Z}\}$ is bounded;
    \item \label{It3}
    There exist $C>0$ and $\gamma\in \R$ such that 
   $\|\Phi_{Z}(t,0)\|\le C e^{\gamma t}$ for every $t\ge 0$ and every $Z\in {\cal S}_{{\cal Z},0}$.
\end{enumerate}
\end{lemma}
\begin{proof}
Let us first prove that Property~\ref{It1} implies Property~\ref{It2}. 
For that, we assume that $\{Z^1_2\cdots Z^k_2\mid k\geq 1,\; (Z^1_1,Z^1_2),\dots, (Z^k_1,Z^k_2)\in {\cal Z}\}$ is unbounded and we are going to prove that $\lam(\Xi_{{\cal Z},0})=+\infty$.
Let $C>0$ be such that
$\|Z_1\|\le C$ for every $(Z_1,Z_2)\in\mathcal{Z}$. Then 
$\|e^{t Z_1}\|\le  e^{C |t|}$ for every $(Z_1,Z_2)\in {\cal Z}$ and $t\in \R$. For every $n\in \mathbb{N}$, let $(Z^{1,n}_1,Z^{1,n}_2),\dots,(Z^{k_n,n}_1,Z^{k_n,n}_2)\in {\cal Z}$ be such that 
\[ \|Z_2^{k_n,n}\cdots Z_2^{1,n}\|\ge  e^{n^2 C}.\]
Then $Z_2^{k_n,n}\cdots Z_2^{1,n}e^{n Z_1^{1,n}}=\Pi_{\omega^n_{k_n}}$ for some $\omega^n\in \Omega_{\mathcal{Z},0}$ with $|\omega^n_{k_n}|=n$. 
Notice that
\begin{align*}
e^{n^2 C}&\le \|Z_2^{k_n,n}\cdots Z_2^{1,n}\|\le  \|\Pi_{\omega^n_{k_n}}\| \|e^{-nZ_1^{1,n}}\|\\
&\le \|\Pi_{\omega^n_{k_n}}\|  e^{nC}.
\end{align*}

Hence
\[\lam(\Xi_{{\cal Z},0})\ge \limsup_{n\to+\infty}\dfrac{\log(\|\Pi_{\omega^n_{k_n}}\|)}{n}=+\infty.\]

Assume now that Property~\ref{It2} holds true and let us prove Property~\ref{It3}. 
Define 
\[v(x)=\sup_{k\ge 0,\;(Z_1^n,Z_2^n)_{n\in\mathbb{N}}\in {\cal Z}^\mathbb{N}}|Z_2^k\cdots Z_2^1x|,\]
with the convention that $Z_2^k\cdots Z_2^1=I_d$ if $k=0$. 
Notice that 
$v$ is finite by Property~\ref{It2} and that
$v(x)\ge |x|>0$ for $x\ne 0$. Moreover, $v$ is  {\color{black}absolutely} homogeneous and satisfies the triangle inequality, hence it is a norm.
Denote by $\|\cdot\|_v$ the matrix norm induced by $v$. Then $\|Z_2\|_v\le 1$ for every $(Z_1,Z_2)\in {\cal Z}$. 
This implies that there exists $\gamma>0$ such that
$\|\Phi_{Z}(t,0)\|_v\le  e^{\gamma t}$ for every $t\ge 0$ and $Z\in {\cal S}_{{\cal Z},0}$. 
Property~\ref{It3} follows. 

The fact that Property~\ref{It3} implies Property~\ref{It1} 
follows from  
inequality \eqref{Xi<=Sigma} with $\tau=0$.
\end{proof}

\begin{remark}
By Remark~\ref{remark-infini} and Lemma~\ref{lem:triip}, $\lambda(\Sigma_{\mathcal{Z},\tau})=+\infty$
if and only if $\tau=0$ and $\{Z^1_2\cdots Z^k_2\mid k\geq 1,\; (Z^1_1,Z^1_2),\dots, (Z^k_1,Z^k_2)\in {\cal Z}\}$ is unbounded.
\end{remark}

{\color{black}
\subsection{A technical result}
}

In the proof of Theorem~\ref{thm:supsup}, 
we will make use of the following technical result, providing a useful bound on the norm of an exponential matrix, which is a variation  of~\cite[Equation (2.11)]{vanloan}. 
\begin{lemma}\label{lemma-est-exp}
If $M\in M_d(\R)$ and $t\geq 0$ then
\[\|e^{tM}\| \leq e^{t\alpha(M)}\sum_{k=0}^{d-1}\frac{t^k d^k \|M\|^k}{k!}. \]
\end{lemma}
\begin{proof}
By Schur triangularization theorem~\cite[Theorem 2.3.1]{horn-johnson}, there exists a unitary matrix $U$ such that we can write
\[M = U^* T U\]
for some upper triangular matrix $T\in M_d(\mathbb{C})$.
 We can write $T=D+N$ where $D$ is the diagonal part of $T$ and $N$ is strictly upper triangular.
As in~\cite[Equation (2.11)]{vanloan} we have
\begin{equation}
\label{vanloan}
\|e^{t M}\|\leq e^{t\alpha(M)}\sum_{k=0}^{d-1}\frac{t^k \|N\|^k}{k!}.
\end{equation}
Moreover, considering the matrix norm $\|A\|_{\infty} = \max_{i,j}|A_{ij}|$, we have
\begin{equation}
\label{norms}
   \|N\|\leq d\|N\|_{\infty}\leq d\|D+N\|_{\infty}\leq d\|D+N\| = d\|M\|, 
\end{equation}
where the first inequality is obtained as a simple application of Cauchy--Schwarz inequality, and the last equality follows from the fact that the transformation $U$ is unitary.
The lemma follows by combining~\eqref{vanloan} with~\eqref{norms}.
\end{proof}

{\color{black}
\subsection{Proof of Theorem~\ref{thm:supsup} and Corollary~\ref{thm:ESEU}}
}
{\color{black}
Theorem~\ref{thm:supsup} is a direct consequence of Remark~\ref{rmk:rho-sup-limsup} together with the last equality in  Proposition~\ref{lambda-rho} below.
}  

\begin{proposition}\label{lambda-rho}
Assume that $\lambda(\Sigma_{\mathcal{Z},\tau})<+\infty$. Then 
\begin{align*}
\lambda(\Sigma_{{\cal Z},\tau})&=\max\left\{\lam(\Xi_{{\cal Z},\tau}),\sup_{(Z_1,Z_2)\in\mathcal{Z}}\alpha(Z_1)\right\}\\
&=\max\left\{\mu(\Xi_{{\cal Z},\tau}),\sup_{(Z_1,Z_2)\in\mathcal{Z}}\alpha(Z_1)\right\}.
\end{align*}
\end{proposition}

\begin{proof}
We first notice that, for every $(Z_1,Z_2)\in \mathcal{Z}$, the flow corresponding to the constant 
signal $Z(\cdot)\equiv (Z_1,Z_2)$, without  switchings, 
satisfies $\Phi_Z(t,0)=e^{t Z_1}$ for every $t\ge 0$. Hence $\alpha(Z_1)\le \lambda(\Sigma_{\mathcal{Z},\tau})$. 
Thus 
\begin{equation}\label{eq:partealpha}
\sup_{(Z_1,Z_2)\in\mathcal{Z}}\alpha(Z_1)\le \lambda(\Sigma_{\mathcal{Z},\tau})
.
\end{equation}
Moreover, since $ \rho(M)\le \|M\|$ for every $M\in M_d(\R)$, and by~\eqref{Xi<=Sigma},
it follows that 
\begin{equation}\label{ineq-oneside}
\mu(\Xi_{{\cal Z},\tau})\le \lam(\Xi_{{\cal Z},\tau})\le  \lambda(\Sigma_{{\cal Z},\tau}).
\end{equation}
Furthermore, by definition, there exist sequences of elements $Z_1^n\in M_d(\R)$, $\tau_n>0$, $\omega^n\in \Omega$, and $k_n\in\mathbb{N}$ with $\lim_{n\to\infty} \left(\tau_n+|\omega^n_{k_n}|\right) = +\infty$ such that
\begin{align*}
\lambda(\Sigma_{\mathcal{Z},\tau})&=\lim_{n\to\infty}\frac{\log(\|e^{\tau_n Z_1^n}\Pi_{\omega_{k_n}^n}\|)}{\tau_n+|\omega^n_{k_n}|}\\&
=\lim_{n\to\infty}\frac{\log(\|e^{\tau_n Z_1^n}\|)+\log(\|\Pi_{\omega_{k_n}^n}\|)}{\tau_n+|\omega^n_{k_n}|}.
\end{align*}
If the sequence $\tau_n$ is bounded then the previous limit is equal to $\limsup_{n\to\infty}\frac{\log(\|\Pi_{\omega_{k_n}^n}\|)}{|\omega^n_{k_n}|}$, hence it is bounded by $\lam(\Xi_{{\cal Z},\tau})$, while if the sequence $|\omega^n_{k_n}|$ is bounded then the limit is equal to $\limsup_{n\to\infty}\alpha(Z_1^n)$ by Lemma~\ref{lemma-est-exp}, and it is bounded by $\sup_{(Z_1,Z_2)\in\mathcal{Z}}\alpha(Z_1)$. If both $\tau_n$ and $|\omega^n_{k_n}|$ tend to infinity (up to a subsequence) then, using the fact that any ratio $\frac{a_1+a_2}{b_1+b_2}$ with $b_1,b_2$ positive is smaller than $\max\left\{\frac{a_1}{b_1},\frac{a_2}{b_2}\right\}$, we still get 
\begin{align*}
&\lambda(\Sigma_{\mathcal{Z},\tau})\leq\limsup_{n\to\infty}\max
\left\{\frac{\log(\|\Pi_{\omega_{k_n}^n}\|)}{|\omega^n_{k_n}|},\frac{\log(\|e^{\tau_n Z_1^n}\|)}{\tau_n}\right\}\\
&=\max
\left\{\limsup_{n\to\infty}\frac{\log(\|\Pi_{\omega_{k_n}^n}\|)}{|\omega^n_{k_n}|},\limsup_{n\to\infty}\frac{\log(\|e^{\tau_n Z_1^n}\|)}{\tau_n}\right\}\\
&\leq \max
\left\{\lam(\Xi_{{\cal Z},\tau}),\sup_{(Z_1,Z_2)\in\mathcal{Z}}\alpha(Z_1)\right\}.
\end{align*}
Combining with~\eqref{eq:partealpha} and~\eqref{ineq-oneside} we get \[\lambda(\Sigma_{\mathcal{Z},\tau})=\max\left\{\lam(\Xi_{{\cal Z},\tau}),\sup_{(Z_1,Z_2)\in\mathcal{Z}}\alpha(Z_1)\right\},\] 
and we are left to prove 
that $\mu(\Xi_{{\cal Z},\tau})$ is equal to $\lam(\Xi_{{\cal Z},\tau})=\lambda(\Sigma_{{\cal Z},\tau})$ whenever
\begin{equation}\label{eq:useff}
\sup_{(Z_1,Z_2)\in\mathcal{Z}}\alpha(Z_1) < \lam(\Xi_{{\cal Z},\tau}).
\end{equation}
 For this purpose, consider a 
flag of subspaces $\{0\}\subsetneq E_1\subsetneq \cdots \subsetneq E_r=\mathbb{R}^d$ 
such that each $E_j$ is invariant for $\{Z_2e^{t Z_1}\mid t\ge \tau,\;(Z_1,Z_2)\in\mathcal{Z}\}$ and $r$ is maximal among all flags with the same property. 
The flag induces a block triangularization 
\[P A P^{-1}=\left(\begin{smallmatrix}
A_{11}&A_{12}&\cdots&A_{1r}\\
0&A_{22}&\cdots
&A_{2r}\\
\vdots&\ddots&\ddots&\vdots\\
0&\dots&0&A_{rr}
\end{smallmatrix}\right),
\]
with $P$ invertible and independent of $A\in \{Z_2e^{t Z_1}\mid t\ge \tau,\;(Z_1,Z_2)\in\mathcal{Z}\}$. 
Up to a linear change of coordinates, we can assume that $P=I_d$.

Let ${\mathcal N}_{i,\tau}:=\{(Z_2e^{t Z_1})_{ii}\mid t\ge \tau,\;(Z_1,Z_2)\in\mathcal{Z}\}$ and consider the corresponding weighted discrete-time switched system $\Xi_i:=\Xi_{\mathcal{N}_{i,\tau}}$. 
Notice that, by maximality of $r$, ${\mathcal N}_{i,\tau}$ is irreducible.  

For $\omega\in \Omega_{\mathcal{Z},\tau}$ and $k\in\mathbb{N}$, 
 the spectrum of $\Pi_{\omega_k}$  is given by the union of the spectra of $(\Pi_{\omega_k})_{11},\dots,(\Pi_{\omega_k})_{rr}$. 
Hence, 
 $\mu(\Xi_{\mathcal{Z},\tau})=\max_{i=1,\dots,r}\mu(\Xi_{i})$.  
 
 Let us conclude the argument by assuming, for now, that  
 \begin{equation}\label{eq:decomp-JH}
     \lam(\Xi_{\mathcal{Z},\tau})=\max_{i=1,\dots,r}\lam(\Xi_{i}).
 \end{equation}

By Corollary~\ref{Berger-Wang}, for  $i=1,\dots,r$ one has 
$\lam(\Xi_{i})=\max\{\hatlam(\Xi_{i}),\mu(\Xi_{i})\}$.  
 Notice that, for every  $(Z_1,Z_2)\in \mathcal{Z}$ and for $i=1,\dots,r$, 
 \[ \|(Z_2e^{t Z_1})_{ii}\|\le \|Z_2e^{t Z_1}\|\le \|Z_2\| \|e^{t Z_1}\|.\]
 Hence, $\hatlam(\Xi_i)\le \sup_{(Z_1,Z_2)\in \mathcal{Z}}\alpha(Z_1)$ for every $i\in \{1,\dots,r\}$.
 
Since $\lambda(\Sigma_{\mathcal{Z},\tau})=\lam(\Xi_{{\cal Z},\tau})$, 
picking $i$ such that $\lam(\Xi_{{\cal Z},\tau})=\lam(\Xi_i)$, we have that 
either $\lambda(\Sigma_{\mathcal{Z},\tau})=\mu(\Xi_{i})\le \mu(\Xi_{\mathcal{Z},\tau})$ or $\lambda(\Sigma_{\mathcal{Z},\tau})=\hatlam(\Xi_{i})\le \sup_{(Z_1,Z_2)\in \mathcal{Z}}\alpha(Z_1)$, proving the desired inequality.

 We are left to prove that, under assumption \eqref{eq:useff},  equality \eqref{eq:decomp-JH} holds true. 
  Notice that $\lam(\Xi_{i})\leq\lam(\Xi_{\mathcal{Z},\tau})$ for $i=1,\dots, r$. The equality is proved by induction on $r$. The case $r=1$ is trivial.

Assume that the equality holds true for some positive integer $r$ and consider 
$\Xi_{\mathcal{Z},\tau}$  with maximal flag $(E_1,\dots, E_{r+1})$ of length $r+1$.
For $\omega\in \Omega_{\mathcal{Z},\tau}$ and $k\in\mathbb{N}$, write 
\[\Pi_{\omega_k}=\begin{pmatrix}(\Pi_{\omega_k})_{11}&(\Pi_{\omega_k})_{1R}\\
0& (\Pi_{\omega_k})_{RR}\end{pmatrix},\]
where $R$ stands for the $r$uple of indices $(2,\dots,r+1)$. 
Applying the induction hypothesis, one deduces that for every $\e>0$ there exists $C_1(\e)>0$ independent of $\omega$ and $k$ such that $\|(\Pi_{\omega_k})_{RR}\|\leq C_1(\e)e^{(\nu+\e)|\omega_k|}$, where $\nu=\max_{i=2,\dots, r+1}\lam(\Xi_{i})$.
 On the other hand, for every $\e>0$ there exists $C_{2}(\e)$ such that
 \[\| (\Pi_{\omega_{j\to k}})_{11}\|
 \leq C_{2}(\e)e^{(\lam(\Xi_1)+\e)(|\omega_k|-|\omega_j|)}.\]

 Notice, moreover, that
 \begin{align*}
    (\Pi_{\omega_k})_{1R}&=\sum_{j=1}^{k} 
    (\Pi_{\omega_{j\to k}})_{11}
(Z_2^{j} e^{\tau_{j} Z_1^{j}})_{1R} (\Pi_{\omega_{j-1}})_{RR}, 
  \end{align*}
and that there exists $C_3(\e)>0$ such that $\|(Z_2^{j} e^{\tau_{j} Z_1^{j}})_{1R}\|\le C_3 e^{\tau_{j}(\e+\sup_{(Z_1,Z_2)\in \mathcal{Z}}\alpha(Z_1))}$.
One deduces that $\lam(\Xi_{\mathcal{Z},\tau})\le \max\{\lam(\Xi_1),\nu,\sup_{(Z_1,Z_2)\in \mathcal{Z}}\alpha(Z_1)\}+\e$. Since $\e$ is arbitrary and we are assuming \eqref{eq:useff}, this concludes the inductive step.
\end{proof}

\begin{remark}
Proposition~\ref{lambda-rho} may fail to hold when $\lambda(\Sigma_{\mathcal{Z},\tau})=+\infty$. 
Consider for example $d=2$, $\tau=0$ and $\mathcal{Z}$ made by the single element $(Z_1,Z_2)$ with $Z_1=(\begin{smallmatrix}0&1\\0&0\end{smallmatrix})$ and $Z_2=(\begin{smallmatrix}1&1\\0&1\end{smallmatrix})$.
Then it is easy to check that 
$\lambda(\Sigma_{\mathcal{Z},\tau})=+\infty$, while $\alpha(Z_1)=0$ and $\mu(\Xi_{\mathcal{Z},0})=0$. 
\end{remark}

{\color{black}{\it Proof of Corollary~\ref{thm:ESEU}.}}
We start by proving Item~\ref{stability}. On the one hand, if $\Sigma_{{\cal Z},\tau}$ is ES then clearly $\lambda(\Sigma_{{\cal Z},\tau})<0$.
On the other hand, $ \lambda(\Sigma_{{\cal Z},\tau})<0$ implies that, 
for every $\gamma\in (\lambda(\Sigma_{{\cal Z},\tau}),0)$,
there exists $T>0$
such that, for every $t> T$ and every $Z\in {\cal S}_{{\cal Z},\tau}$, $\|\Phi_Z(t,0)\|\le e^{\gamma t}$. We are left to show that $\{\Phi_Z(t,0)\mid Z\in {\cal S}_{{\cal Z},\tau},\;t\in [0,T]\}$ is bounded. 
In the case $\tau>0$ this is straightforward (see Remark~\ref{remark-infini}). If $\tau=0$,
since $ \lam(\Xi_{{\cal Z},0})\le   \lambda(\Sigma_{{\cal Z},0})<+\infty$, we deduce the boundedness  from Property~\ref{It3} of Lemma~\ref{lem:triip}.  

Concerning Item~\ref{instability}, if $\Sigma_{{\cal Z},\tau}$ is EU then clearly $\lambda(\Sigma_{{\cal Z},\tau})>0$. On the other hand, if $\lambda(\Sigma_{{\cal Z},\tau})>0$ then, from Theorem~\ref{thm:supsup}, there exist either $(Z_1,Z_2)\in {\cal Z}$ such that $\alpha(Z_1)>0$  or 
$\omega\in\Omega_{\cal Z,\tau}$ and $k\in\mathbb{N}$ such that $\rho(\Pi_{\omega_k})>1$. 
In the first case, $\Sigma_{{\cal Z},\tau}$ is obviously EU by taking the signal constantly equal to $(Z_1,Z_2)$. In the second case, we consider the sequence $\omega^\star\in \Omega_{\mathcal{Z},\tau}$ obtained by repeating $\omega_k$ infinitely many times. 
Let $Z^\star\in \mathcal{S}_{{\cal Z},\tau}$ be the signal associated with $\omega^\star$ and $T=|\omega_k|$. There exist $c>1$ and $x_0\in \R^d\setminus\{0\}$ such that $|\Phi_{Z^\star}(nT,0)x_0|\ge c^n|x_0|$.
Let $C$ and $\gamma$ be as in Property~3 of Lemma~\ref{lem:triip}. Then 
for $t\in [nT,(n+1)T)$ one has 
\begin{align*}
&|\Phi_{Z^\star}(t,0)x_0|\ge \\
&|\Phi_{Z^\star}((n+1)T,0)x_0|\|\Phi_{Z^\star}((n+1)T,s)^{-1}\|\ge \\
& c^{n+1}|x_0|\frac{e^{-|\gamma| T}}{C}.
\end{align*}
This concludes the proof that 
$\Sigma_{{\cal Z},\tau}$ is EU.  \hfill $\Box$

{\color{black}
\subsection{The case where the Lyapunov exponent is attained on unswitched trajectories}
}
It follows from Proposition~\ref{lambda-rho} that when  $\lam(\Xi_{{\cal Z},\tau})<\lambda(\Sigma_{{\cal Z},\tau})<+\infty$, {\color{black}i.e., when the 
exponential growth rate of $\Sigma_{{\cal Z},\tau}$ is not captured by trajectories that switch infinitely many times, }
 there must exist $(Z_1,Z_2)\in {\cal Z}$ such that $\lam(\Xi_{{\cal Z},\tau})<\alpha(Z_1)$. The following proposition investigates such a situation.

\begin{proposition}
Assume
$\lam(\Xi_{{\cal Z},\tau})<\alpha(Z_1)$ for some $(Z_1,Z_2)\in {\cal Z}$. 
Then $Z_2x=0$ for every generalized eigenvector $x$
of $Z_1$ associated with  an eigenvalue of real part
$\alpha(Z_1)$. 
\end{proposition}
\begin{proof}
For simplicity of notation, we write the proof 
when $x$ is a generalized eigenvector of $Z_1$ associated with a real eigenvalue, the general case being similar.
Recall that $x$ is a generalized eigenvector 
of $Z_1$ associated with the eigenvalue $\alpha(Z_1)$ if there exist $k\geq 1$ linearly independent vectors $x_1,\dots, x_k$ (referred to as a Jordan chain of length $k$) so that $x=x_k$ and  $Z_1x_j=\alpha(Z_1)x_j+\sum_{i=1}^{j-1}x_i$ for $1\leq j\leq k$.

We will prove the conclusion by induction on $k$. For $k=1$, $x$ is simply an eigenvector 
of $Z_1$ associated with the eigenvalue $\alpha(Z_1)$
and we can assume without loss of generality that $|x|=1$. Suppose by contradiction that 
$Z_2x\ne 0$.
Using the fact that $e^{tZ_1}x=e^{\alpha(Z_1)t}x$ and based on the definition of $\lam(\Xi_{{\cal Z},\tau})$, one has that
\begin{align*}
    \alpha(Z_1)&>\lam(\Xi_{{\cal Z},\tau})\ge \limsup_{t\to +\infty} \frac{\log|Z_2e^{tZ_1}x|}{t}\\
    &=\alpha(Z_1)+\limsup_{t\to +\infty} \frac{\log|Z_2x|}{t}=\alpha(Z_1),
\end{align*}
yielding a contradiction. Therefore,  we must have $Z_2x=0$.

Assume now that the conclusion holds true for every $j$ with $1\leq j\leq k-1$. Consider a generalized eigenvector $x$ with Jordan chain $x_1,\dots,x_k$ of length $k$. Notice that for every $j$ with $1\leq j\leq k-1$, $x_j$ is a generalized eigenvector with Jordan chain $x_1,\dots,x_j$ of length $j\leq k-1$. Applying the induction hypothesis one gets that $Z_2x_j=0$. Moreover, for $t\geq 0$,
\[
e^{tZ_1}x=e^{\alpha(Z_1)t}\Big(x+\sum_{j=1}^{k-1}\frac{t^{j-1}}{(j-1)!}x_j\Big).
\]
Then $Z_2e^{tZ_1}x=e^{\alpha(Z_1)t}Z_2x$. If $Z_2x\neq 0$,
we argue as before and reach a contradiction. This ends the induction argument and concludes the proof of the proposition.
\end{proof}

{\color{black}
\subsection{Application of the Berger--Wang formula: continuity of the maximal Lyapunov exponent}\label{ssec:continuity}
}
A consequence of {\color{black}Theorem~\ref{thm:supsup}} is the continuity of the maximal Lyapunov exponent, as detailed in the following proposition.

\begin{proposition}\label{lambda-continuity}
Let $\tau>0$ and ${\cal U}$ be the set of bounded subsets ${\cal Z}$ of $M_d(\R)\times M_d(\R)$.
Endow ${\cal U}$ with the topology induced by the Hausdorff distance. Then 
the function ${\cal Z}\mapsto \lambda(\Sigma_{{\cal Z},\tau})$ is continuous on  ${\cal U}$.
\end{proposition}

{\color{black}
Before proving the proposition, let us provide a statement of converse Lyapunov 
theorem for impulsive switched systems. Such a result can be obtained by  classical arguments (as, for instance, those of \cite[Lemma~6.2]{cai2007smooth}).

\begin{proposition}\label{converse-cont}
System $\Sigma_{{\cal Z},\tau}$ 
satisfies $\lambda(\Sigma_{{\cal Z},\tau})<0$
 if and only if
 $\sup_{(Z_1,Z_2)\in\mathcal{Z}}\alpha(Z_1)<0$ and
there exist $c\ge 1$, $\gamma>0$, and  $V:\R^{d}\to \R_+$ convex, absolutely homogeneous, and Lipschitz continuous such that, for every $x\in \R^{d}$, $(Z_1,Z_2)\in{\cal Z}$ and $t\in  \R_{\ge \tau}$, we have 
\begin{align}
|x|\leq V(x)\leq c|x|, 
 \label{comp-cont}\\
V(Z_2e^{t Z_1}x)\leq e^{-\gamma t}V(x). 
\label{decay-cont}
\end{align}
\end{proposition}

}

{

{\color{black}
{\it Proof of Proposition~\ref{lambda-continuity}.}}
We begin by noting that the map ${\cal Z}\mapsto \sup_{Z\in {\cal Z}}\alpha(Z_1)$ is continuous on ${\cal U}$. This follows directly from the uniform continuity of $\alpha(\cdot)$ on bounded subsets of $M_{d}(\R)$.

We now claim that, if $\lambda(\Sigma_{{\cal Z},\tau})<0$ for some ${\cal Z}\in {\cal U}$, then for $\e>0$ small enough and every ${\cal Z}'\in {\cal U}$ with $d_{H}({\cal Z},{\cal Z}')<\e$ we have $\lambda(\Sigma_{{\cal Z}',\tau})<0$. Indeed,  
{\color{black}
let $\lambda(\Sigma_{{\cal Z},\tau})<0$ and $V$ be as in Proposition~\ref{converse-cont}.}
Let ${\cal Z}'$ be sufficiently close to ${\cal Z}$ such that $\sup_{Z'\in {\cal Z}'}\alpha(Z_1')<\frac12\sup_{Z\in {\cal Z}}\alpha(Z_1)<0$.
Then, by Lemma~\ref{lemma-est-exp},  there exist $C\ge 1$ and $\tilde\gamma\in (0,\gamma)$ such that for every $t\geq 0$, $Z\in {\cal Z}$, and $Z'\in {\cal Z}'$, one has 
\begin{equation}\label{eq:estZ0}
\|e^{tZ_1}\|\leq Ce^{-\tilde\gamma t}\qquad
\textrm{ and} \qquad 
\|e^{tZ'_1}\|\leq Ce^{-\tilde\gamma t}.
\end{equation}
Let us now show that, for every $Z'\in {\cal Z}'$ and every $t\ge \tau$, there exists $Z\in {\cal Z}$ such that 
\begin{equation}\label{eq:estZ}
e^{-\gamma t}+L\|Z'_2e^{tZ'_1}-Z_2e^{tZ_1}\|\le e^{-\frac{\tilde\gamma}{2} t},
\end{equation}
{\color{black}where $L>0$ is such that $V$ is $L$-Lipschitz continuous.}
To see that, 
notice that
for every $Z\in {\cal Z}$, $Z'\in {\cal Z}'$, and $t\geq \tau$, we have 
\begin{align}\label{eq:inter0}
&e^{-(\gamma-\frac{\tilde\gamma}2)t}+
e^{\frac{\tilde\gamma}{2} t}
L\|Z'_2e^{tZ'_1}-Z_2e^{tZ_1}\|\le \\
&e^{-\frac{\tilde\gamma}{2}\tau}+
e^{\frac{\tilde\gamma}2t}
\|Z'_2(e^{tZ'_1}-e^{tZ_1})\|
+\|Z'_2-Z_2\|\|e^{\frac{\tilde\gamma}2t}e^{tZ_1}\|.\nonumber
\end{align}
Using \eqref{eq:estZ0}, note that, for every $t\geq 0$, one has that 
\[
e^{\frac{\tilde\gamma}2t}
\|Z'_2(e^{tZ'_1}-e^{tZ_1})\|\leq C_1e^{-\frac{\tilde\gamma}2t}
\; \textrm{ and }\;
\|e^{\frac{\tilde\gamma}2t}e^{tZ_1}\|\leq 
C.
\]
Pick $\kappa\in (0,\frac{1-e^{-\frac{\tilde\gamma}{2}\tau}}2)$. 
Fix $T>0$ so that $C_1e^{-\frac{\tilde\gamma}2t}\leq \kappa$ for all $t\ge T$ 
and choose $Z$ close enough to $Z'$ so that 
\[\sup_{t\in[0,T]}\|Z'_2(e^{tZ'_1}-e^{tZ_1})\|\leq \kappa.
\]
Hence we deduce that for every $t\geq \tau$, $e^{\frac{\tilde\gamma}2t}
\|Z'_2(e^{tZ'_1}-e^{tZ_1})\|\leq \kappa$. Similarly, choose again $Z$ close enough to $Z'$ so that $C\|Z'_2-Z_2\|\leq \kappa$. Collecting all the above estimates, one gets that 
the left-hand side of \eqref{eq:inter0} is upper bounded by $e^{-\frac{\tilde\gamma}{2}\tau}+2\kappa$ for $t\geq \tau$, which implies \eqref{eq:estZ}.
{\color{black}Using the $L$-Lipschitz continuity of $V$, it follows from~\eqref{decay-cont} that, for every $x\in \R^{d}$,  $(Z'_1,Z'_2)\in{\cal Z}'$, and $t\ge\tau$, one has
\begin{align*} 
V(Z'_2e^{tZ'_1}x)&\le V(Z_2e^{tZ_1}x)+L\|Z'_2e^{tZ'_1}x-Z_2e^{tZ_1}x\|\\
&\leq e^{-\gamma t}V(x)+L\|Z'_2e^{tZ'_1}x-Z_2e^{tZ_1}\|V(x)\\
&\le e^{-\frac{\tilde\gamma}{2} t}V(x),
\end{align*}
where we have used the fact that $\|x\|\le V(x)$.}
{\color{black}
Proposition~\ref{converse-cont}
implies that }
$\lambda(\Sigma_{{\cal Z}',\tau})<0$. 
This concludes the proof of the claim.

The above claim actually proves the 
lower semi-continuity of ${\cal Z}\mapsto \lambda(\Sigma_{{\cal Z},\tau})$. 
Indeed, let $\delta>0$ and fix ${\cal Z}\in {\cal U}$.  Define $\xi=-\lambda(\Sigma_{{\cal Z},\tau})-\delta$. Then, by Lemma~\ref{shift}, we have $\lambda(\Sigma_{{\cal Z}^{\xi},\tau})=-\delta<0$.  The claim guarantees that there exists a neighbourhood ${\cal W}$ of ${\cal Z}$ in ${\cal U}$ 
such that for all ${\cal Z}'\in {\cal W}$ we have $\lambda(\Sigma_{{{\cal Z}'}^{\xi},\tau})<0$, which means that 
$\lambda(\Sigma_{{\cal Z}',\tau})<\lambda(\Sigma_{{\cal Z},\tau})+\delta$.

For the upper semi-continuity of ${\cal Z}\mapsto \lambda(\Sigma_{{\cal Z},\tau})$, consider $\delta>0$ and ${\cal Z}\in {\cal U}$.  According to Theorem~\ref{thm:supsup}, we have
\begin{align*}
    \lambda(\Sigma_{{\cal Z},\tau})=&\max\Big\{\sup_{(Z_1,Z_2)\in {\cal Z}}\alpha(Z_1),\\
&\sup_{Z\in {\cal S}_{{\cal Z},\tau},\;k\in \Theta^\star(Z)}\frac{\log(\rho(\Phi_Z(t_k,0)))}{t_k}\Big\}. 
\end{align*}
If $\lambda(\Sigma_{{\cal Z},\tau})=\sup_{(Z_1,Z_2)\in {\cal Z}}\alpha(Z_1)$, then, for every ${\cal Z}'$ close enough to ${\cal Z}$ in ${\cal U}$  we have 
\begin{align*}
\lambda(\Sigma_{{\cal Z}',\tau})&\ge\sup_{(Z'_1,Z'_2)\in {\cal Z}'}\alpha(Z'_1)\geq \sup_{(Z_1,Z_2)\in {\cal Z}}\alpha(Z_1)-\delta\\
&=\lambda(\Sigma_{{\cal Z},\tau})-\delta.
\end{align*} 
Otherwise, there exist $Z\in {\cal S}_{{\cal Z},\tau}$ and $k\in \Theta^\star(Z)$ such that 

\[\frac{\log(\rho(\Phi_Z(t_k,0)))}{t_k}>\lambda(\Sigma_{{\cal Z},\tau})-\frac{\delta}{2}.\]

Let ${\cal W}$ be a neighbourhood of ${\cal Z}$  in ${\cal U}$ such that for every ${\cal Z}'\in {\cal W}$ there exists $Z'\in {\cal S}_{{\cal Z}',\tau}$  such that $\Theta^\star(Z')=\Theta^\star(Z)$ and 
\[\frac{\log(\rho(\Phi_{Z'}(t_k,0)))}{t_k}
\ge \frac{\log(\rho(\Phi_{Z}(t_k,0)))}{t_k}-\frac{\delta}{2}.\]
It follows that 
\begin{align*}
&\lambda(\Sigma_{{\cal Z}',\tau})\ge
\sup_{\tilde Z\in {\cal S}_{{\cal Z}',\tau},\;\tilde k\in \Theta^\star(\tilde Z)}
\frac{\log(\rho(\Phi_{\tilde Z}(t_{\tilde k},0)))}{t_{\tilde k}}\\
&\ge 
\frac{\log(\rho(\Phi_{Z'}(t_k,0)))}{t_k}
\ge \frac{\log(\rho(\Phi_{Z}(t_k,0)))}{t_k}-\frac{\delta}{2}.
\end{align*}
It follows that $\lambda(\Sigma_{{\cal Z}',\tau})\ge \lambda(\Sigma_{{\cal Z},\tau})-\delta$. In both cases, we conclude that there exists a neighbourhood ${\cal W}$ of ${\cal Z}$ in ${\cal U}$ such that $\lambda(\Sigma_{{\cal Z}',\tau})>\lambda(\Sigma_{{\cal Z},\tau})-\delta$ for every ${\cal Z}'\in {\cal W}$.
\hfill $\Box$
}

{

If $\tau=0$ 
we can extend Proposition~\ref{lambda-continuity} as follows.

\begin{remark}
Let
$\tau=0$ and
$\lambda(\Sigma_{{\cal Z},0})<+\infty$.
Set ${\cal Z}_2=\{Z_2\mid (Z_1,Z_2)\in {\cal Z}\}$ and let 
${\cal U}$ be the set of bounded subsets ${\cal Z}'$ of $M_d(\R)\times M_d(\R)$ such that
\[\{Z_2'\mid (Z_1',Z_2')\in {\cal Z}'\}={\cal Z}_2.\]
Notice that, by Lemma~\ref{lem:triip}, $\lambda(\Sigma_{{\cal Z}',0})<+\infty$ for every ${\cal Z}'\in {\cal U}$. 
Then  showing the continuity of ${\cal U}\ni {\cal Z}'\mapsto \lambda({\cal Z}')$ follows the same lines as that of Proposition~\ref{lambda-continuity}, where the argument to get \eqref{eq:estZ} in the case $\tau=0$ now requires the extra fact that for every $T>0$ and every bounded subset ${\cal B}$ of $M_d(\R)$ there exists $c>0$ such that  $\|e^{t Z_1}-e^{t Z_1'}\|\le c t\|Z_1-Z_1'\|$ for every $Z_1,Z_1'\in {\cal B}$ and $t\in[0,T]$.

\end{remark}
}

\section{Conclusion}
This paper addresses the stability analysis of impulsive linear switched systems by examining their representation as weighted discrete-time switched systems. {\color{black}We provide 
a Berger--Wang-type result that establishes the equality between two measures of asymptotic growth, one based on the operator norm and the other on the spectral radius. 
This result is} first developed in the general context of weighted discrete-time switched systems and subsequently applied to impulsive linear switched systems.  In particular, using the Berger--Wang formula, we establish a characterization of exponential stability for impulsive linear switched systems in terms of the sign of their maximal Lyapunov exponent.

\bibliographystyle{plain}        
\bibliography{biblio}

\end{document}